\documentclass[english]{amsart}
\usepackage[T1]{fontenc}
\usepackage[latin1]{inputenc}
\usepackage{amssymb}

\makeatletter
 \theoremstyle{plain}    
 \newtheorem{thm}{Theorem}[section]
 \numberwithin{equation}{section} %% Comment out for sequentially-numbered
 \numberwithin{figure}{section} %% Comment out for sequentially-numbered
 \theoremstyle{plain}
 \theoremstyle{remark}    
 \newtheorem{notation}[thm]{Notation} 
 \theoremstyle{plain}    
 \newtheorem{prop}[thm]{Proposition} %%Delete [thm] to re-start numbering
 \theoremstyle{definition}
 \newtheorem{defn}[thm]{Definition}
 \theoremstyle{remark}
 \newtheorem{rem}[thm]{Remark}
 \theoremstyle{plain}    
 \newtheorem{lem}[thm]{Lemma} %%Delete [thm] to re-start numbering
 \theoremstyle{plain}    
 \newtheorem{cor}[thm]{Corollary} %%Delete [thm] to re-start numbering
 \theoremstyle{remark}    
 \newtheorem*{acknowledgement*}{Acknowledgement} 

\usepackage{babel}
\makeatother
\begin{document}

\title{An Asymptotics of Variance of the Lattice Points Count}

\author{Ji\v r\'{i} Jan\'a\v cek}

\address{Institute of Physiology, Academy of Sciences of the Czech Republic,
V\'{i}de\v nsk\'a 1083, 142 20 Praha, Czech Republic}

\begin{abstract}
The variance of the number of lattice points inside the dilated bounded
set $rD$ with random position in $\mathbf{\mathbb{R}}^{d}$ has asymptotics
$\sim r^{d-1}$ if the rotational quadratic average of the modulus
of the Fourier transform of the set is $O\left(\rho^{-d-1}\right)$.
The asymptotics follows from a Wiener's Tauberian theorem.
\end{abstract}



\maketitle

\section{Introduction}

The number of lattice points in a set with random position can be
used for estimation of the volume of the set and its variance has
been studied for a long time \cite{5,6}. If $\mathbf{T}$ is a $d$-periodical
point lattice of spatial intensity $\alpha$ in the $d$-dimensional
Euclidean space $\mathbf{\mathbb{R}}^{d}$ then the mean value of\[
\left(\mathbf{card}\left(\left(B_{d}\left(r\right)+x\right)\cap\mathbf{T}\right)-\alpha\lambda^{d}\left(B_{d}\left(r\right)\right)\right)^{2},\]
where $\lambda^{d}$ is Lebesgue measure, i.e. the variance of the
lattice points count in the ball of radius $r$ with uniform random
position, is

\[
C_{\mathbf{T}}H^{d-1}\left(\partial B_{d}\left(r\right)\right)\Phi\left(r\right).\]
Here $C_{\mathbf{T}}$ is a lattice constant, $H^{d-1}$ is the surface
measure and $\Phi$ defined by the above equality fulfills $\lim_{t\rightarrow\infty}\frac{1}{t}\int_{0}^{t}\Phi\left(t\right)\: dt=1$
and is bounded, $\limsup_{t\rightarrow\infty}\Phi\left(t\right)\leq2$
\cite{3,4}. Hence the variance of the lattice point count in the
ball has asymptotics {}``in the mean'' $C_{\mathbf{T}}H^{d-1}\left(\partial B_{d}\left(1\right)\right)r^{d-1}$
and is $O\left(r^{d-1}\right)$, $r\rightarrow\infty$.

Let $D$ be a compact body the isotropic covariogram $\overline{\gamma}_{D}$
(rotational average of $\gamma_{D}=I_{D}\star I_{-D}$) of which is
a fractional integral of the Weyl type of a finite signed measure
$\sigma$ on $\mathbf{\mathbb{R}}^{+}$: \[
\overline{\gamma}_{D}\left(t\right)=\frac{1}{\Gamma\left(\frac{d+1}{2}+1\right)}\int_{t}^{\infty}\left(s-t\right)^{\frac{d+1}{2}}d\sigma\left(s\right)\]
for $t\in\mathbf{\mathbb{R}}^{+}$. This means that the fractional
derivative of $\overline{\gamma}_{D}$ of the order $\frac{d+1}{2}$
has bounded variation. Then the asymptotics {}``in the mean'' of
the form $C_{\mathbf{T}}H^{d-1}\left(\partial D\right)r^{d-1}$ of
the variance of the lattice points count inside randomly rotated and
shifted set $rD$ follows easilly from Theorem 2.8 in \cite{key-13}.

The variance of the lattice points count inside randomly rotated and
shifted bounded set $D$ dilated by $r>0$ is $O\left(r^{d-1}\right)$,
$r\rightarrow\infty$, if the rotational quadratic average of the
modulus of the Fourier transform of the set $D$, defined as $\int_{S_{d-1}}\left|\widehat{I_{D}}\left(\rho\xi\right)\right|^{2}d\xi$
is $O\left(\rho^{-d-1}\right)$, $\rho\rightarrow\infty$ \cite{key-12};
the later property was proved for convex sets and sets with $C^{\frac{3}{2}}$
boundary \cite{key-6}. The aim of this paper is to prove that this
assumption on the growth of the modulus of the Fourier transform also
yields the asymptotics of the form $C_{\mathbf{T}}H^{d-1}\left(\partial D\right)r^{d-1}$
of the variance of the lattice points count inside randomly placed
bounded full $d$-dimensional set $rD$ with sufficiently regular
boundary (locally finite union of sets of finite reach \cite{8}).
The above asymptotics can be thus used for compact bodies with piecewise
$\mathbf{C}^{2}$ smooth boundary and $d$-dimensional convex sets
in $\mathbf{\mathbb{R}}^{d}$.

\section{The Variance of the Lattice Points Count}

\begin{notation}
\label{d:periodic-meas}Let $\mathbf{T}$ be a $d$-periodic lattice
of points in the $d$-dimensional Euclidean space $\mathbf{\mathbb{R}}^{d}$
defined by the regular matrix $A\in\mathbf{\mathbb{R}}^{d\times d}$
as $\mathbf{T}\left(A\right)=A\mathbf{\mathbb{Z}}^{d}$, where $\mathbf{\mathbb{Z}}^{d}$
is set of all points in $\mathbf{\mathbb{R}}^{d}$ with integral co-ordinates.
$\mathbf{T}$ has the fundamental region $F_{\mathbf{T}}=A\left[0,1\right)^{d}$
of volume $\lambda^{d}\left(F_{\mathbf{T}}\right)=\mathrm{det}\: A$,
where $\lambda^{d}$ is the Lebesgue measure; hence the spatial intensity
of $\mathbf{T}$ is $\alpha=\left(\mathrm{det}\: A\right)^{-1}$.
The group dual to the group $\mathbf{T}\left(A\right)$ is $\mathbf{T}^{*}=A^{-1}\mathbf{\mathbb{Z}}^{d}$.
\end{notation}
Fourier transform of a function $f\in\mathbf{L}^{1}\left(\mathbf{\mathbb{R}}^{d}\right)$
is\begin{equation}
\widehat{f}\left(\xi\right)=\int_{\mathbf{\mathbb{R}}^{d}}f\left(x\right)\exp\left(-2\pi ix\cdot\xi\right)dx.\label{eq:Four-transf}\end{equation}

If $f$ is moreover spherically symmetric then $r^{d-1}f\left(r\right)\in\mathbf{L}^{1}\left(\mathbf{\mathbb{R}}^{+}\right)$
and Fourier transform of $f$ can be expressed as the Haenkel transform

\emph{\begin{equation}
\widehat{f}\left(\rho\right)=2\pi\rho^{1-\frac{d}{2}}\int_{0}^{\infty}r^{\frac{d}{2}}J_{\frac{d}{2}-1}\left(2\pi\rho r\right)f\left(r\right)dr,\label{eq:Haenkel}\end{equation}
}where $J_{\frac{d}{2}-1}$ is the Bessel function of the first kind.

$\kappa_{d}=\pi^{\frac{d}{2}}\Gamma\left(\frac{d}{2}+1\right)^{-1}$
is the volume of the unit ball $B_{d}\left(1\right)$ in $\mathbf{\mathbb{R}}^{d}$,
where $\Gamma$ is the Euler gamma function.

$I_{D}$ is the characteristic function of the set $D$.

\begin{prop}
\label{t:p-grid-var}Let $\mathbf{T}$ be a $d$-periodic lattice
of points and let $D$ be a bounded measurable set in $\mathbf{\mathbb{R}}^{d}$.
Then \begin{equation}
\int_{F_{\mathbf{T}}}\left(\mathbf{card}\left(\left(D+x\right)\cap\mathbf{T}\right)\right)\alpha\: dx=\alpha\lambda^{d}\left(D\right),\label{eq:mean}\end{equation}
and\begin{equation}
\begin{array}{c}
\int_{F_{\mathbf{T}}}\left(\mathbf{card}\left(\left(D+x\right)\cap\mathbf{T}\right)-\alpha\lambda^{d}\left(D\right)\right)^{2}\alpha\: dx=\\
={\displaystyle \sum_{0\neq\xi\in\mathbf{T}^{*}}}\left|\widehat{I_{D}}\left(\xi\right)\right|^{2},\end{array}\label{eq:var-freq}\end{equation}
where $\alpha$ is the spatial intensity of $\mathbf{T}$.
\end{prop}
\begin{proof}
Equation (\ref{eq:mean}) can be proved by standard arguments and
equation (\ref{eq:var-freq}) follows from the Parseval theorem, see
Theorem 2.3 in \cite{key-13}.
\end{proof}
\begin{defn}
Covariogram of a bounded measurable set $D$ is the function \[
\gamma_{D}\left(x\right)=I_{D}\star I_{-D}\left(x\right)=\int_{\mathbf{\mathbb{R}}^{d}}I_{D}\left(x\right)I_{D}\left(x-y\right)dy.\]
 It follows from the properties of Fourier transform that $\widehat{\gamma_{D}}=\left|\widehat{I_{D}}\right|^{2}$
is a nonnegative function. The isotropic covariogram is $\overline{\gamma_{D}}\left(\left|x\right|\right)=\int_{\mathbf{SO}_{d}}\gamma_{MD}\left(x\right)dM$,
where $MD$ is the set $D$ rotated by $M\in\mathbf{S}\mathbf{O}_{d}$
and the integration uses the invariant probabilistic measure on $\mathbf{S}\mathbf{O}_{d}$,
the group of rotations in $\mathbf{\mathbb{R}}^{d}$; an equivalent
definition is $\overline{\gamma_{D}}\left(u\right)=\int_{S_{d-1}}\gamma_{D}\left(ux\right)dx$.
The Haenkel transform of the isotropic covariogram is $\widehat{\overline{\gamma_{D}}}$.
\end{defn}
\begin{rem}
\label{rem:Loewe} As was already discussed in Remark 2.5 in \cite{key-13},
it follows from this definition that $\gamma_{D}$ is bounded and,
as $\widehat{\gamma_{D}}\geq0$, the function $\widehat{\gamma_{D}}$
is integrable in $\mathbf{\mathbb{R}}^{d}$ (see \cite{key-1} Theorem
9.). Further, $\rho^{d-1}\widehat{\overline{\gamma_{D}}}\left(\rho\right)\geq0$
is integrable in $\mathbf{\mathbb{R}}^{+}$ by Fubini theorem. $\gamma_{D}$
is then the inverse Fourier transform (\ref{eq:Four-transf}) of $\widehat{\gamma_{D}}$
(\cite{key-1} Theorem 8.) and $\overline{\gamma_{D}}$ is the (inverse)
Haenkel transform (\ref{eq:Haenkel}) of $\widehat{\overline{\gamma_{D}}}\left(\rho\right)$.
We have from (\ref{eq:var-freq}) and (\ref{eq:mean}) using the variance
decomposition lemma \cite{7} \[
\begin{array}{c}
\int_{\mathbf{SO}_{d}}\int_{F_{\mathbf{T}}}\left(\mathbf{card}\left(\left(MD+x\right)\cap\mathbf{T}\right)-\alpha\lambda^{d}\left(D\right)\right)^{2}\alpha\: dx\: dM=\\
={\displaystyle \sum_{0\neq\xi\in\mathbf{T}^{*}}}\widehat{\overline{\gamma_{D}}}\left(\left|\xi\right|\right).\end{array}\]

\end{rem}
\begin{lem}
\label{lemma:Wiener-Tauber}Let $r^{2}f\left(r\right)\geq0$ be bounded
measurable function on $\mathbf{\mathbb{R}}^{+}$. Then \[
\lim_{R\rightarrow\infty}\frac{1}{R}\int_{0}^{R}r^{2}f\left(r\right)\: dr=\lim_{h\rightarrow0+}\frac{\pi^{-\frac{3}{2}}}{h}\frac{\Gamma\left(\frac{d+1}{2}\right)}{\Gamma\left(\frac{d}{2}\right)}\int_{0}^{\infty}\left(1-\Gamma\left(\frac{d}{2}\right)\frac{J_{\frac{d}{2}-1}\left(2\pi hr\right)}{\left(\pi hr\right)^{\frac{d}{2}-1}}\right)f\left(r\right)\: dr\]
whenever at least one of the two limits exists.
\end{lem}
\textbf{Proof.} The case $d=1$ was proved in \cite{key-4} (Theorem
21) using Wiener's Tauberian theorem (see e.g. \cite{key-3}): Let
$\varphi\in\mathbf{L}^{\infty}\left(\mathbf{\mathbb{R}}\right)$,
$K\in\mathbf{L}^{1}\left(\mathbf{\mathbb{R}}\right)$, $\widehat{K}$
has no root and $K\star\varphi\left(t\right)\rightarrow a\widehat{K}\left(0\right)$
as $t\rightarrow\infty$. Then for each $g\in\mathbf{L}^{1}\left(\mathbf{\mathbb{R}}\right)$
$g\star\varphi\left(t\right)\rightarrow a\widehat{g}\left(0\right)$
as $t\rightarrow\infty$. What follows is an extension of this proof
to higher dimensions.

By the substitution $r=$$\exp(t)$ and defining

\[
\begin{array}{cc}
R=\exp\left(\eta\right)=\frac{1}{h},\; & r^{2}f\left(r\right)=\varphi\left(t\right),\end{array}\]
we obtain equivalent formulation of the theorem suitable for direct
application of the Wiener's Tauberian theorem

\[
\lim_{\eta\rightarrow\infty}\int_{-\infty}^{\infty}K_{1}\left(\eta-t\right)\varphi\left(t\right)\: dt=\lim_{\eta\rightarrow\infty}\int_{-\infty}^{\infty}K_{2}\left(\eta-t\right)\varphi\left(t\right)\: dt\]
whenever at least one of the two limits exists. Here \[
\begin{array}{cc}
K_{1}\left(s\right)=I_{\left\{ s|s>0\right\} }\exp\left(-s\right),\; & K_{2}\left(s\right)=\exp\left(s\right)L\left(\exp\left(-s\right)\right)\end{array},\]
 \[
L\left(u\right)=\pi^{-\frac{3}{2}}\frac{\Gamma\left(\frac{d+1}{2}\right)}{\Gamma\left(\frac{d}{2}\right)}\left(1-\Gamma\left(\frac{d}{2}\right)\frac{J_{\frac{d}{2}-1}\left(2\pi u\right)}{\left(\pi u\right)^{\frac{d}{2}-1}}\right).\]
We can easilly see that $\widehat{K_{1}}\left(\tau\right)=\int_{0}^{\infty}\exp\left(-s-2\pi is\tau\right)ds=\frac{1}{1+2\pi i\tau}$
has no real root and $\widehat{K_{1}}\left(0\right)=1$. If moreover
$\widehat{K_{2}}\left(0\right)=\int_{0}^{\infty}u^{-2}L\left(u\right)du=1$
and the left side limit exists, then also the right side limit exists
and the equality follows from the Wiener's Tauberian theorem. To prove
the statement in the oposite direction it remains to show that $\widehat{K_{2}}\left(\tau\right)=\int_{0}^{\infty}u^{-2+2\pi i\tau}L\left(u\right)du$
has no real root. 

To complete the proof by establishing the validity of assumptions
concerning the function $\widehat{K_{2}}$, we will evaluate the integral
\[
\int_{0}^{\infty}u^{-2+2\pi i\tau}\left(1-\Gamma\left(\frac{d}{2}\right)\frac{J_{\frac{d}{2}-1}\left(2\pi u\right)}{\left(\pi u\right)^{\frac{d}{2}-1}}\right)du.\]
Using $\int t^{-\nu}J_{\nu+1}\left(t\right)dt=-t^{-\nu}J_{\nu}\left(t\right)$
we get

\[
2\pi\Gamma\left(\frac{d}{2}\right)\int_{0}^{\infty}u^{-2+2\pi i\tau}\int_{0}^{u}\left(\pi s\right)^{1-\frac{d}{2}}J_{\frac{d}{2}}\left(2\pi s\right)ds\: du\]
and integration by parts gives

\[
\frac{2\pi\Gamma\left(\frac{d}{2}\right)}{1-2\pi i\tau}\left(-\left[u^{-1+2\pi i\tau}\int_{0}^{u}\frac{J_{\frac{d}{2}}\left(2\pi s\right)}{\left(\pi s\right)^{\frac{d}{2}-1}}ds\right]_{0}^{\infty}+\pi^{1-2\pi i\tau}\int_{0}^{\infty}\frac{J_{\frac{d}{2}}\left(2\pi u\right)}{\left(\pi u\right)^{\frac{d}{2}-2\pi i\tau}}du\right).\]

The first term is zero, as the limits in $0$ and infinity are zero
by l'Hospital formula and asymptotic properties of the Bessel function.
From the formula $\int_{0}^{\infty}t^{a}J_{\nu}\left(2t\right)dt=\frac{1}{2}\Gamma\left(\frac{\nu+a+1}{2}\right)\Gamma\left(\frac{\nu-a+1}{2}\right)^{-1}$,
valid if $Re\: a<\frac{1}{2}$, $Re\: a+\nu>-1$, it follows that

\[
\widehat{K_{2}}\left(\tau\right)=\frac{\pi^{-2\pi i\tau-\frac{1}{2}}}{1-2\pi i\tau}\frac{\Gamma\left(\frac{d+1}{2}\right)\Gamma\left(\frac{1}{2}+\pi i\tau\right)}{\Gamma\left(\frac{d+1}{2}-\pi i\tau\right)}.\]

Now it is easy to see that $\widehat{K_{2}}\left(0\right)=1$, and
$\widehat{K_{2}}$ has no real root because gamma function has none
and because all poles of gamma function are negative.

\begin{rem}
\label{rem:Rataj} If $D$ is a bounded full-dimensional locally finite
union of sets of finite reach, then the derivative of the covariance
from the right $\overline{\gamma_{D}}'^{+}\left(0\right)=-\frac{\kappa_{d-1}}{d\kappa_{d}}H^{d-1}\left(\partial D\right)$
\cite{8}.
\end{rem}
\begin{thm}
\label{t:asymp-var}Let $\mathbf{T}$ be a $d$-periodic lattice of
points, $r\in\mathbf{\mathbb{R}}^{+}$, $D$ a bounded measurable
set such that finite $\overline{\gamma_{D}}'^{+}\left(0\right)$ exists
and $\Phi$ a function on $\mathbf{\mathbb{R}}^{+}$ defined by equation\[
\int_{\mathbf{SO}_{d}}\int_{F_{\mathbf{T}}}\left(\mathbf{card}\left(\left(rMD+x\right)\cap\mathbf{T}\right)-\alpha\lambda^{d}\left(D\right)\right)^{2}\alpha\: d\lambda^{d}\left(x\right)dM=\]
\[
=\frac{-\overline{\gamma_{D}}'^{+}\left(0\right)}{2\pi^{2}\kappa_{d-1}}\left(\sum_{0\neq\xi\in\mathbf{T}^{*}}\left|\xi\right|^{-d-1}\right)\Phi\left(r\right)r^{d-1}.\]
Further, let $\rho^{d+1}\widehat{\overline{\gamma_{D}}}\left(\rho\right)$
be bounded measurable on $\mathbf{\mathbb{R}}^{+}$. Then $\Phi$
is bounded and\begin{equation}
\lim_{t\rightarrow\infty}\frac{1}{t}\int_{0}^{t}\Phi\left(t\right)\: dt=1.\label{eq:int-Phi}\end{equation}
If moreover $\lim_{\rho\rightarrow\infty}\rho^{d+1}\widehat{\overline{\gamma_{D}}}\left(\rho\right)$
exists then 

\begin{equation}
\lim_{t\rightarrow\infty}\Phi\left(t\right)=1.\label{eq:lim-Phi}\end{equation}

\end{thm}
\begin{proof}
From Remark \ref{rem:Loewe} we have \[
\int_{\mathbf{SO}_{d}}\int_{F_{\mathbf{T}}}\left(\mathbf{card}\left(\left(rMD+x\right)\cap\mathbf{T}\right)-\alpha\lambda^{d}\left(D\right)\right)^{2}\alpha\: d\lambda^{d}\left(x\right)dM=\sum_{0\neq\xi\in\mathbf{T}^{*}}\widehat{\overline{\gamma_{D}}}\left(r\left|\xi\right|\right)r^{2d}.\]

We shall prove first that the auxiliary function $\Psi$ defined by
the equation\[
-\overline{\gamma_{D}}'^{+}\left(0\right)\Psi\left(t\right)=2\pi^{2}\kappa_{d-1}t^{d+1}\widehat{\overline{\gamma_{D}}}\left(t\right)\]
 has the property (\ref{eq:int-Phi}) or (\ref{eq:lim-Phi}). It is
easy to see that the function\[
\Phi\left(t\right)=\frac{{\displaystyle \sum_{0\neq\xi\in\mathbf{T}^{*}}}\left|\xi\right|^{-d-1}\Psi\left(t\left|\xi\right|\right)}{{\displaystyle \sum_{0\neq\xi\in\mathbf{T}^{*}}}\left|\xi\right|^{-d-1}}\]
is bounded and has then the same property (\ref{eq:int-Phi}) or (\ref{eq:lim-Phi})
as the auxiliary function $\Psi$.

Let $\rho^{d+1}\widehat{\overline{\gamma_{D}}}\left(\rho\right)$
be bounded measurable on $\mathbf{\mathbb{R}}^{+}$. From the Lemma
\ref{lemma:Wiener-Tauber}, Proposition \ref{eq:Haenkel} and Remark
\ref{rem:Loewe} it follows that \[
\lim_{R\rightarrow\infty}\frac{1}{R}\int_{0}^{R}2\pi^{2}\kappa_{d-1}\rho^{d+1}\widehat{\overline{\gamma_{D}}}\left(\rho\right)\: d\rho=\]
 \[
=\lim_{h\rightarrow0+}\frac{1}{h}\int_{0}^{\infty}\left(d\kappa_{d-1}-2\pi\left(h\rho\right)^{1-\frac{d}{2}}J_{\frac{d}{2}-1}\left(2\pi h\rho\right)\right)\rho^{d-1}\widehat{\overline{\gamma_{D}}}\left(\rho\right)\: d\rho=\]
\[
=\lim_{h\rightarrow0+}\frac{1}{h}\left(\overline{\gamma_{D}}\left(0\right)-\overline{\gamma_{D}}\left(h\right)\right)=-\overline{\gamma_{D}}'^{+}\left(0\right)\]
and if the left side limit exists then

\[
\lim_{r\rightarrow\infty}2\pi^{2}\kappa_{d-1}r^{d+1}\widehat{\overline{\gamma_{D}}}\left(r\right)=-\overline{\gamma_{D}}'^{+}\left(0\right).\]
 
\end{proof}
\begin{cor}
\label{n:var_coef_reg} From Theorem \ref{t:asymp-var} and Remark
\ref{rem:Rataj} it follows that if $D$ is a bounded full-dimensional
locally finite union of sets of finite reach such that $\rho^{d+1}\widehat{\overline{\gamma_{K}}}\left(\rho\right)$
is bounded (or has a limit in $+\infty$), then\begin{equation}
\begin{array}{c}
\int_{\mathbf{SO}_{d}}\int_{F_{\mathbf{T}}}\left(\mathbf{card}\left(\left(rMD+x\right)\cap\mathbf{T}\right)-\alpha\lambda^{d}\left(D\right)\right)^{2}\alpha\: d\lambda^{d}\left(x\right)dM=\\
=C_{\mathbf{T}}H^{d-1}\left(\partial D\right)\Phi\left(r\right)r^{d-1},\end{array}\label{eq:asymp}\end{equation}
where

\[
C_{\mathbf{T}}=\frac{1}{2\pi^{2}d\kappa_{d}}\;\sum_{0\neq n\in\mathbf{\mathbb{Z}^{d}}}\left|A^{-1}n\right|^{-d-1}\]
is a lattice constant and $\Phi$ fulfills (\ref{eq:int-Phi}) (or
has limit $1$).
\end{cor}

\section{Discussion}

The assumption of the Theorem \ref{t:asymp-var} and Corollary \ref{n:var_coef_reg},
that $\widehat{\overline{\gamma_{D}}}\left(\rho\right)$ is $O\left(\rho^{-d-1}\right)$,
$\rho\rightarrow\infty$, holds for convex sets and sets with $C^{\frac{3}{2}}$
boundary (the boundary of the set can be decomposed into finitely
many neighbourhoods such that given any pair of points $P,Q$ in the
neighbourhood, $\left|\left(P-Q\right)n\left(Q\right)\right|\leq c\left|P-Q\right|^{\frac{3}{2}}$,
where $n\left(Q\right)$ is a unit normal to the set in $Q$) \cite{key-6}.
The asymptotics $\sim r^{d-1}$ of the variance of the lattice points
count in the mean value (\ref{eq:asymp}) thus holds for sets with
piecewise $\mathbf{C}^{2}$ smooth boundary and convex sets in $\mathbf{\mathbb{R}}^{d}$.

As was already said in the introduction, Theorem \ref{t:asymp-var}
gives results similar to Theorem 2.8 in \cite{key-13} for a compact
body $D$ with smooth isotropic covariogram $\overline{\gamma}_{D}$.
Slightly different situation is studied in \cite{key-13}: locally
finite periodic measure is studied instead of the point lattice and
the size of the body is fixed while the scale $s$ of the lattice
tends to zero and the counting measure is multiplied by factor $s^{d}$;
consequently, the asymptotics $\sim s^{-d-1}$ is obtained there.
As in \cite{key-13} we may generalize the results of the present
paper to locally finite periodic measures. In \cite{key-13}, lattice
constants $C_{\mathbf{T}}$ in Equation \ref{eq:asymp} were calculated
for some important point lattices.

\begin{acknowledgement*}
This study was supported by the Academy of Sciences of the Czech Republic,
grant No. A100110502 and AV0Z 50110509.
\end{acknowledgement*}

\end{document}